\documentclass[12pt,leqno]{article}
\usepackage{amsmath, amsthm, amsfonts, amssymb, color}
\usepackage{mathrsfs}
\usepackage[notcite,notref]{showkeys}

\theoremstyle{definition}

\newcommand{\scr}[1]{\mathscr #1}
\definecolor{wco}{rgb}{0.5,0.2,0.3}

\numberwithin{equation}{section}

\newcommand{\ua}{\uparrow}

\title{{\bf Pointwise  Characterizations  of   Curvature and  Second Fundamental Form on Riemannian Manifolds}\footnote{Supported in
 part by  NNSFC (11771326, 11431014).} }
\author{
{\bf     Feng-Yu Wang$^{a,c)}$  and Bo Wu$^{b)}$ }\\
\footnotesize{  $^{a)}$Center of Applied Mathematics, Tianjin University, Tianjin 300072, China}\\
\footnotesize {$^{b)}$ School  of Mathematical Sciences, Fudan
University, Shanghai 200433, China}\\
 \footnotesize{ $^{c)}$Department of Mathematics,
Swansea University, Singleton Park, SA2 8PP, United Kingdom}\\
\footnotesize{  wangfy@bnu.edu.cn, F.-Y.Wang@swansea.ac.uk, wubo@fudan.edu.cn}}

\date{}

\begin{document}

\maketitle

\def\R{\mathbb R} \def\EE{\mathbb E} \def\Z{\mathbb Z} \def\ff{\frac} \def\ss{\sqrt}
\def\H{\mathbb H}
\def\dd{\delta} \def\DD{\Delta} \def\vv{\varepsilon} \def\rr{\rho}
\def\<{\langle} \def\>{\rangle} \def\GG{\Gamma} \def\gg{\gamma}
\def\ll{\lambda} \def\LL{\Lambda} \def\nn{\nabla} \def\pp{\partial}
\def\d{\text{\rm{d}}} \def\Id{\text{\rm{Id}}}\def\loc{\text{\rm{loc}}} \def\bb{\beta} \def\aa{\alpha} \def\D{\scr D}
\def\E{\scr E} \def\si{\sigma} \def\ess{\text{\rm{ess}}}
\def\beg{\begin} \def\beq{\beg}  \def\F{\scr F}
\def\Ric{\text{\rm{Ric}}}
\def\Var{\text{\rm{Var}}}
\def\Ent{\text{\rm{Ent}}}
\def\Hess{\text{\rm{Hess}}}\def\B{\scr B}
\def\e{\text{\rm{e}}} \def\ua{\underline a} \def\OO{\Omega} \def\b{\mathbf b}
\def\oo{\omega}     \def\tt{\tilde} \def\Ric{\text{\rm{Ric}}}
\def\cut{\text{\rm{cut}}} \def\P{\mathbb P} \def\ifn{I_n(f^{\bigotimes n})}
\def\fff{f(x_1)\dots f(x_n)} \def\ifm{I_m(g^{\bigotimes m})} \def\ee{\varepsilon}
\def\C{\scr C}
\def\M{\scr M}\def\ll{\lambda}
\def\X{\scr X}
\def\T{\scr T}
\def\A{\mathbf A}
\def\LL{\scr L}\def\LLL{\Lambda}
\def\gap{\mathbf{gap}}
\def\div{\text{\rm div}}
\def\dist{\text{\rm dist}}
\def\cut{\text{\rm cut}}
\def\supp{\text{\rm supp}}
\def\Cov{\text{\rm Cov}}
\def\Dom{\text{\rm Dom}}
\def\Cap{\text{\rm Cap}}\def\II{{\mathbb I}}\def\beq{\beg{equation}}
\def\sect{\text{\rm sect}}\def\H{\mathbb H}

\begin{abstract} Let $M$ be  a complete Riemannian manifold  possibly with a boundary $\pp M$. For any $C^1$-vector field $Z$, by  using gradient/functional  inequalities of the (reflecting) diffusion process  generated by $L:=\DD+Z$,   pointwise characterizations are presented for the Bakry-Emery curvature of  $L$ and the second fundamental form of  $\pp M$ if exists.  These extend and strengthen the recent results derived by A. Naber for the uniform norm $\|\Ric_Z\|_\infty$ on manifolds without boundary. A key point of the present study is to apply the asymptotic formulas for these two tensors found by the first named author, such that the proofs are significantly simplified.
\end{abstract}

\vskip 2cm

\section{Introduction}\label{sect1}

Let $M$ be a $d$-dimensional complete  Riemannian manifold possibly with a boundary
$\partial M$. Let $L=\DD+Z$ for a $C^1$ vector field $Z$.  We intend to characterize  the Bakry-Emery curvature $\Ric_Z:= \Ric-\nn Z$  and the second fundamental form $\II$ of the boundary $\pp M$ using the (reflecting) diffusion process generated by $L$. When $\pp M=\emptyset$, we set $\II=0.$

There are many equivalent characterizations for the (pointwise or uniform) lower bound of $\Ric_Z$ and $\II$ using gradient/functional inequalities of the (Neumann) semigroup generated by $L$, see e.g. \cite{W2} and references within. However, the corresponding upper bound characterizations are still open.
It is   known that for stochastic analysis on the path space, one needs conditions on the norm of $\Ric_Z$, see \cite{CHL, CW, D, F, Hsu0,W1,W11}  and references within.  Recently,   A. Naber \cite{N,HN} proved that the uniform bounded condition on $\Ric_Z$ for $Z=-\nabla f$ is equivalent to some gradient/functional inequalities on the path space, and thus  clarified the necessity of bounded conditions used in the above mentioned references.  In this paper, we aim to present pointwise characterizations for the norm of $\Ric_Z$ and $\II$ when $\pp M\ne \emptyset,$ which  allow these quantities unbounded on the manifold.

Let $(X_t^x)_{t\ge 0}$  be the (reflecting if $\pp M$ exists) diffusion process generated by $L=\DD+Z$ on $M$ starting at point $x$, and let $(U_t^x)_{t\ge 0}$ be the horizontal lift onto the frame bundle $O(M):= \cup_{x\in M} O_x(M)$, where $O_x(M)$ is the set of all orthonormal basis of the tangent space $T_x M$ at point $x$. It is well known that $(X_t^x, U_t^x)_{t\ge 0}$ can be constructed as the unique solution to the  SDEs:
\begin{equation}\label{eq1.1}\beg{split} & \d X_t^x= \ss 2\, U_t^x\circ d W_t+ Z(X_t^x)\d t + N(X_t^x)\d l_t^x,\ \ X_0^x=x,\\
&\d U_t^x=\sqrt{2}\, H_{U_t^x}(U_t^x)\circ\d W_t+H_{Z}(U_t^x)\d t+H_N(U_t^x)\d l_t^x,\ \ U_0^x\in O_x(M),\end{split}
\end{equation}
where $W_t$ is the $d$-dimensional Brownian motion on a complete filtration probability space $(\OO, \{\F_t\}_{t\ge 0}, \P)$, $N$ is the inward unit normal vector field of $\pp M$, $H_\cdot: TM\to TO(M)$ is the horizontal lift, $H_{u}:=(H_{ue_i})_{1\le i\le d}$ for $u\in O(M)$ and the canonical orthonormal basis
$\{e_i\}_{1\le i\le d}$ on $\R^d$,  and $l_t$ is an adapted increasing process which increases only when $X_t^x\in \pp M$ which is   called the local time of $X_t^x$ on $\pp M$.  In the first part of this paper, we assume that the solution is non-explosive, so that the (Neumann) semigroup $P_t$ generated by $L$  is given by
$$P_tf(x)= \EE f(X_t^x),\ \ x\in M, f\in \B_b(M), t\ge 0.$$

For a fixed $T>0$, consider the path space
 $W_T(M):=C([0,T];M) $
and the class of smooth cylindric   functions
$$\aligned \F C^\infty_{T}:=\Big\{F(\gamma)&=f(\gamma_{t_1},\cdots,\gamma_{t_m}):\ m\geq1,~\gamma\in W_T(M),\\&
~~~~~~~~~~~~~~~~~~~~~~
0<t_1<t_2\cdots<t_m\leq T, ~f\in C_0^\infty(M^m)\Big\}.\endaligned$$  Let  $$\H_T=
 \left\{h\in C([0,T];\mathbb{R}^d): \ \ h(0)=0,~
\|h\|^2_{\mathbb{H}_T}:=\int_0^T|h_s'|^2\d s<\infty\right\}.$$ For any  $F\in\F C^\infty_T$ with
$F(\gamma)=f\big(\gamma(t_1),\cdots,\gamma(t_m)\big)$, the Malliavin gradient $DF(X_{[0,T]}^x)$ is an $\H_T$-valued random variable satisfying
\beq\label{MG}\beg{split} &\dot{D}_sF(X_{[0,T]}^x):=\ff{\d}{\d s} DF(X_{[0,T]}^x)\\
&= \sum_{t_i> s}  (U_{t_i}^x)^{-1}\nn_if\big(X_{t_1}^x,\cdots,X_{t_m}^x\big),\ \ s\in [0,T],\end{split} \end{equation}where $\nabla_i$ is the (distributional) gradient operator for the $i$-th
component on $M^m$, and $P_u:\R^d\to \R^d$ is the projection along $u^{-1}N$, i.e.
$$\<P_ua, b\>:=\<ua,N\>\<ub,N\>,\quad  a, b\in\R^d, u\in \cup_{x\in \pp M}O_x(M).$$ Note that 

For $K\in C(M; [0,\infty))$ and $\si\in C(\pp M; [0,\infty))$, we introduce  the following random measure $\mu_{x,T}$ on $[0,T]$:
\beq\label{MU}
  \mu_{x,T}(\d s):= \e^{\int^s_0K(X_r^x)\d r+\int^s_0\sigma(X_r^x)\d l_r^x}\big\{K(X_s^x) \d s+  \sigma(X_s^x) \d l_s^x\big\}. \end{equation} For any $t\in [0,T]$, consider  the energy form
$$\E^{K,\sigma}_{t,T}(F,F)
 = \EE\bigg\{\big(1+\mu_{x,T}([t,T])\big)
  \bigg(\big|\dot{D}_tF(X_{[0,T]}^x)\big|^2
 +\int^{T}_{t}   \big|\dot{D}_sF(X_{[0,T]}^x)\big|^2\mu_{x,T}(\d s)\bigg)\bigg\}$$ for $ F\in \scr F C_{T}^\infty.$
Our main result is the following.

\beg{thm}\label{T1.1}  Let $K\in C(M; [0,\infty))$ and $\si\in C(\pp M; [0,\infty))$ be such that
\beg{align} \label{DE} \EE\e^{(2+\vv)\int_0^T \{K(X_s^x)\d s +\si(X_s^x)\d l_s^x\}}<\infty \  {\rm for\ some}\ \vv,T>0.\end{align} For any $p, q\in [1,2]$,   the following statements are equivalent each other:
\beg{enumerate}\item[$(1)$] For any $x\in M$ and $y\in\pp M$,
\beg{align*} &\|\Ric_Z\|(x):= \sup_{X\in T_x M, |X|=1} \big|\Ric(X,X)-\<\nn_XZ, X\>\big|(x)\le K(x),\\
&\|\II\|(y):= \sup_{Y\in T_y\pp M, |Y|=1}|\II(Y,Y)|(y) \le \si(y).\end{align*}
\item[$(2)$]  For any $f\in C_0^\infty(M)$, $T>0$, and $x\in M$,
\beg{align*} &|\nn P_T f|^p(x) \le \EE\Big[(1+\mu_{x,T}([0,T]))^p|\nn f|^p(X_T^x)\Big],\\
&\Big|\nn  f(x)- \ff 1 2   \nn P_Tf(x) \Big|^q
 \le \EE\bigg[\big(1+\mu_{x,T}([0,T])\big)^{q-1} \\
 &\qquad \times \bigg(\Big|\nn f(x)-\ff 1 2 U_0^x(U_T^x)^{-1} \nn f(X_T^x)\Big|^q
  +\ff {\mu_{x,T}([0,T])} {2^q}  \big| \nn f(X_T^x)\big|^q  \bigg)\bigg].\end{align*}
\item[$(3)$] For any $F\in \scr F C_T^\infty, x\in M$ and $T>0$,
\beg{align*} \big|\nn_x \EE F(X_{[0,T]}^x)\big|^q \le \EE&\bigg[\big(1+\mu_{x,T}([0,T])\big)^{q-1}\\
&\ \times \bigg(\big|\dot D_0 F(X_{[0,T]}^x)\big|^q
  + \int_0^T  \big|\dot D_s F(X_{[0,T]}^x)\big|^q\mu_{x,T}(\d s)\bigg)\bigg].\end{align*}
\item[$(4)$] For any  $t_0,t_1\in [0,T]$ with $t_1>t_0$, and any $ x\in M$,  the following log-Sobolev inequality holds:
$$\aligned &\EE\left[\EE\big(F^2(X_{[0,T]}^x)|\F_{t_1}\big) \log \EE(F^2(X_{[0,T]}^x)|\F_{t_1})\right]\\&-\EE\left[\EE\big(F^2(X_{[0,T]}^x)|\F_{t_0}\big) \log \EE(F^2(X_{[0,T]}^x)|\F_{t_0})\right]
 \le 4 \int_{t_0}^{t_1} \E_{s,T}^{K,\si}(F,F)\d s,\  \  F\in \scr F C_T^\infty.\endaligned$$
 \item[$(5)$] For any  $t \in [0,T]$  and $ x\in M$,  the following Poincar\'e inequality holds:
$$   \EE \Big[\big\{\EE(F(X_{[0,T]}^x)|\F_{t})\big\}^2\Big]-   \Big\{\EE  \big[F (X_{[0,T]}^x)\big] \Big\}^2
 \le 2\int_0^t \E_{s,T}^{K,\si}(F, F)\d s,\ \ F\in \scr F C_T^\infty.$$   \end{enumerate}\end{thm}

\paragraph{Remark 1.1.} (1) When $\pp M=\emptyset, Z=-\nn f$ and $K$ is a constant, it is proved in \cite[Theorem 2.1]{N} that $\|\Ric_Z\|_\infty\le K$ is equivalent to each of (3)-(5) with $\si=0$ and a slightly different formulation of $\E_{s,T}^{K,0}$. Comparing with these equivalent statements using references functions on the path space, the statement (2) only depends on   reference functions on  $M$ and    is thus easier to verify.

(2) An important problem in geometry is to identify the Ricci curvature, for instance, to characterize Einstein manifolds where $\Ric$ is a  constant tensor. According to Theorem \ref{T1.1}, $\Ric$ is identified by $\nn Z$ if and only if all/some of items (2)-(5) hold for $K=0$.

\

We will prove this result in the next section. In Section 3, the equivalence of $(1)$,    (4) and (5) are proved without condition \eqref{DE} but using  the class of truncated cylindrical functions replacing $\F C_T^\infty$.

\section{Proof}

 We first introduce some known results from the monograph  \cite{W2} which hold  under a  condition weaker than   \eqref{DE}.

Let $f\in C_0^\infty(M)$ with $|\nn f(x)|=1$ and  $\Hess_f(x)=0$.  According to \cite[Theorem 3.2.3]{W2},  if $x\in M\setminus \pp M$ then for any $p>0$ we have
\beg{equation}\beg{split} \label{RIC}   \Ric_Z(\nn f,\nn f)(x)  &= \lim_{t\downarrow 0} \ff{P_t|\nn f|^p(x)-|\nn P_t f|^p(x)}{pt} \\
&= \lim_{t\downarrow 0} \ff 1 t \bigg(\ff{P_tf^2(x)-(P_tf)^2(x)}{2t}-|\nn P_tf(x)|^2\bigg);\end{split}\end{equation} and by \cite[Theorem 3.2.3]{W2},
if $x\in \pp M$ and $\nn f\in T_x \pp M$  then
\beg{equation}\label{II} \beg{split} \II(\nn f,\nn f)(x)&=\lim_{t\downarrow 0}  \ff{\ss\pi}{2p\ss t} \Big\{P_t|\nn f|^p(x)-|\nn P_t f|^p(x)\Big\} \\
&=\lim_{t\downarrow 0} \ff{3\ss\pi}{8\ss t}\bigg(\ff{P_tf^2(x)-(P_tf)^2(x)}{2t}-|\nn P_tf|^2(x)\bigg).\end{split}\end{equation}
We note that in \cite[(3.2.9)]{W2}, $\ss\pi$ is misprinted as $\pi$.

Next, let $\Ric_Z(u)$ for $u\in O(M)$ and $\II(u), P_u$ for $u\in \cup_{x\in \pp M} O_xM$ are matrix-valued functions with
\beg{align*} &\<P_ua, b\>=\<ua,N\>\<ub,N\>,\\
&  \<\Ric_Z(u)a,b\>:=\Ric_Z(ua, ub),\\
& \<\II(u)a,b\>:=\II\big(ua-\<ua,N\>N,  ub-\<ub,N\>N\big),\ \ a,b\in \R^d.\end{align*}
According to \cite[Lemma 4.2.3]{W2},
for any $F\in \scr FC_T^\infty$ with $F(\gg)= f(\gg_{t_1},\cdot, \gg_{t_N}), f\in C_0^\infty(M)$ and $0\le t_1<\cdots \leq t_N$,
\beq\label{GR} (U_0^x)^{-1} \nn_x  \EE\big[F(X_{[0,T]}^x)\big] = \sum_{i=1}^N \EE\big[Q_{0,t_i}^x(U_{t_i}^x)^{-1} \nn_i f(X_{t_1}^x,\cdots,X_{t_N}^x)\big],\end{equation}
where $\nn_x$ denotes the gradient in $x\in M$ and $\nn_i$ is the gradient with respect to the $i$-th component, and for any $s\ge 0$, $(Q_{s,t}^x)_{t\ge s}$ is an adapted right-continuous process on $\R^d\otimes \R^d$ satisfies $Q_{s,t}^x P_{U_t^x}=0$ if $X_t^x\in\pp M$ and
\beq\label{QQ} Q_{s,t}^x= \bigg(I-\int_s^t Q_{s,r}^x\big\{\Ric_Z(U_r^x)\d r + \II(U_r^x)\d l_r^x\big\}\bigg)
\Big(I-1_{\{X_t^x\in \pp M\}}P_{U_t^x}\Big). \end{equation} The multiplicative functional $Q_{s,t}^x$ was introduced by Hsu \cite{H1} to investigate gradient estimate on $P_t$. For convenience, let $Q_t^x:=Q_{0,t}^x$. In particular, taking $F(\gg)= f(\gg_t)$ in \eqref{GR},  we obtain
\beq\label{GR2} \nn P_t f(x)= U_0^x \EE\big[Q_t^x(U_t^x)^{-1}\nn f(X_t^x)\big],\ \ x\in M, f\in C_0^\infty(M), t\ge 0.\end{equation}

Finally, for the above  $F\in \scr F C_T^\infty$, let
\beq\label{TTD} \tt D_tF(X_{[0,T]}^x)= \sum_{i: t_i> t}Q_{t,t_i}^xU_{t_i}^{-1} \nn_if(X_{t_1}^x,\cdots, X_{t_N}^x),\ \
t\in [0,T].\end{equation}
Then    \cite[Lemma  4.3.2]{W2} (see also \cite{W11}) implies that
\beq\label{MF} \EE\big(F(X_{[0,T]}^x)\big|\F_t\big)= \EE[F(X_{[0,T]}^x)] +\ss 2 \int_0^t\Big\<\EE(\tt D_sF(X_{[0,T]}^x) |\F_s), \d W_s\Big\>,\ \ t\in [0,T].\end{equation}

\beg{proof}[Proof of Theorem \ref{T1.1}] It is well known that the log-Sobolev inequality in (4) implies the Poincar\'e inequality in (5), below we prove the theorem by verifying the following implications respectively:
  (1) $\Rightarrow$ (3) for all $q\ge 1$;   (3) $\Rightarrow$ (2) for all $p=q$;
 (2) for some $p\ge 1$ and $q\in [1,2] \Rightarrow$ (1);   (5) $\Rightarrow$ (1); and
  $(1)\Rightarrow (4).$

For simplicity, below we will 
write $F$ and $ f$ for  $F(X_{[0,T]}^x)$ and $f(X_{t_1}^x,\cdots, X_{t_N}^x)$ respectively.

  {\bf (a)}  (1) $\Rightarrow$ (3) for all $q\ge 1$. By \eqref{MG}, \eqref{GR} and \eqref{QQ} we have
\beg{equation*} \beg{split}
&U_0^{-1}\nabla_x\mathbb{E}[F]=\mathbb{E}\bigg[\sum_{i=1}^NQ^{x}_{t_i}(U_{t_i}^x)^{-1}\nabla_if\bigg]
 \\&
=\mathbb{E}\bigg[\sum_{i=1}^N\Big(I-\int^{t_i}_0Q^x_s\Ric_Z(U_s)\d s
-\int^{t_i}_0Q^x_s\mathbb{I}_{U^x_s}\d l_s^x\Big) (U_{t_i}^x)^{-1}\nabla_if\bigg]\\
&=\mathbb{E}\bigg[\sum_{i=1}^N (U_{t_i}^x)^{-1}\nabla_if\\&~~~-\sum_{i=1}^N\Big(\int^{t_i}_0Q^x_s\Ric_Z(U_s^x)\d s
+\int^{t_i}_0Q^x_s\mathbb{I}_{U^x_s}\d l_s^x\Big) (U_{t_i}^x)^{-1}\nabla_if\bigg]\\
\\&=\mathbb{E}\bigg[\dot{D}_0F-\int^T_0\big\{Q^x_s\Ric_Z(U_s^x)\dot{D}_sF\big\}\d s-\int^T_0\big\{Q^x_s\mathbb{I}(U^x_s)\dot{D}_sF\big\}\d l_s^x\bigg].\end{split} \end{equation*}
By   \cite[Theorem 3.2.1]{W2}, we have
\beq\label{QQ2} \left\|Q^x_s\right\|\leq \exp\left[\int^s_0K(X_r)\d r+\int^s_0\sigma(X_r)\d l_r^x\right].
\end{equation}  Combining these with (1), \eqref{MU}, and using H\"older's inequality twice, we obtain
\begin{align*}
&\big|\nabla_x\mathbb{E}[F]\big|^q \leq \bigg\{\EE |\dot{D}_0F|
+\EE \int^T_0 |\dot{D}_sF|\mu_{x,T}(\d s)\bigg\}^q\\
&\le  \EE\bigg\{ |\dot{D}_0F| + \int_0^T
  |\dot{D}_sF| \mu_{x,T}(\d s)\bigg\}^q\\
  &\le \EE \bigg\{\bigg(|\dot{D}_0F|^q +\ff {\big(\int_0^T
  |\dot{D}_sF(X_{[0,T]}^x)|\mu_{x,T}(\d s) \big)^q} {\{\mu_{x,T}([0,T]) \}^{q-1}}  \bigg)\big(1+\mu_{x,T}([0,T])\big)^{q-1}\bigg\}\\
  &\le \EE \bigg\{\bigg( |\dot{D}_0F|^q +\int_0^T
  |\dot{D}_sF(X_{[0,T]}^x)|^q \mu_{x,T}(\d s)\bigg)\big(1+\mu_{x,T}([0,T])\big)^{q-1}\bigg\}.\end{align*}
Thus,   the inequality in (3) holds.

{\bf (b)} (3) $\Rightarrow$ (2) for all $p=q$. Take $F(\gg)= f(\gg_T)$. Then $\EE F(X_{[0,T]}^x)=P_T f(x)$ and by \eqref{MG},   $|\dot D_s F|\le |\nn f(X_T)|$ for $s\in [0,T].$ So,
the first inequality  in (2) with $p=q$ follows from (3) immediately. Similarly, by taking $F(\gg)=f(\gg_0)-\ff 1 2 f(\gg_T)$, we have
 $ \EE F= f(x) -\ff 1 2 P_T f(x)$ and \beg{align*}
 & |\dot D_0 F|= \Big|\nn f(x)-\ff 1 2 U_0^x (U_T^x)^{-1}\nn f(X_T^x)\Big|,\\
  &|\dot D_s F|\le \ff 1 2 |\nn f(X_T^x)|,\ \ s\in (0,T].\end{align*}  Then   the second inequality in (2) is implied by (3).

{\bf (c)} (2) for some $p\ge 1$ and $q\in [1,2] \Rightarrow$ (1). Let $x\in M\setminus \pp M$. There exists $r>0$ such that $B(x,r):=\{y\in M: \rr(x,y)\le r\}\subset M\setminus \pp M$, where $\rr$ is the Riemannian distance. Let $\tau_r=\inf\{t\ge 0: \rr(x,X_t^x)\ge r\}$. By \cite[Lemma 3.1.1]{W2} (see also \cite[Lemma 2.3]{ATW09}), there exists a constant $c>0$ such that
\beq\label{LO} \P(\tau_r\le T) \le \e^{-c/T},\ \ T\in (0,1].\end{equation}  Then $\P(l_T^x>0)\le  \e^{-c/T}$ so that for each $n\ge 1$
\beq\label{LI} \lim_{T\to 0} T^{-n} \l_T^x=0,\quad \P-a.s..\end{equation} Combining this with \eqref{MU} we obtain
\beq\label{MUL} \lim_{T\to 0} \ff{\mu_{x,T}([0,T])}T= K(x).\end{equation} Therefore, by the dominated convergence theorem due to \eqref{DE}, the first inequality in (2) and \eqref{RIC} yield
\beq\label{X1} \beg{split}  -\Ric_Z(\nn f,\nn f)(x)&= \lim_{T\to 0}\ff{|\nn P_Tf|^p(x)-P_T|\nn f|^p(x)}{pT}\\
 &\le  \lim_{T\to 0}\ff{\EE\big\{[(1+\mu_{x,T}([0,T]))^p-1]|\nn f|^p(X_T^x)\big\}}{pT} = K(x),\end{split}\end{equation}  where $f\in C_0^\infty(M)$ with $\Hess_f(x)=0$ and $|\nn f(x)|=1$.
This implies  $\Ric_Z(X,X)\ge -K(x)$  for any $X\in T_x M$ with $|X|=1$.

Next, we prove that the second inequality in (2) implies $\Ric_Z\le K$. By H\"older's inequality,   the second inequality in (2) for some $q\in [1,2]$ implies the same inequality  for $q=2$:
\beg{align*} &\Big|\nn f(x)-\ff 1 2 \nn P_Tf(x)\Big|^2\\
&\le \EE\left[(1+\mu_{x,T}([0,T]))\Big(\Big|\nn f(x)-\ff 1 2 U_0^x (U_T^x)^{-1} \nn f(X_T^x)\Big|^2 +\ff{\mu_{x,T}([0,T])} 4 |\nn f(X_T^x)|^2\Big)\right].\end{align*} Then
\beq\label{X2} \beg{split} &\ff{|\nn P_T f(x)|^2-P_T|\nn f(x)|^2}{4T}
 \le \ff 1 T\EE\bigg\{\big\<\nn f(x), \nn P_T f(x) -\EE [U_0^x (U_T^x)^{-1} \nn f(X_T^x)\big\> \\
&\qquad \qquad\qquad\qquad\qquad\qquad~~+\mu_{x,T}([0,T])\Big|\nn f(x)-\ff 1 2 U_0^x (U_T^x)^{-1}\nn f(X_T^x)\Big|^2 \\
&\qquad \qquad\qquad\qquad\qquad\qquad~~+\ff{(1+\mu_{x,T}([0,T]))\mu_{x,T}([0,T])} 4 |\nn f(X_T^x)|^2\bigg\}.\end{split}\end{equation}
Combining this with \eqref{RIC} and \eqref{MUL}, we arrive at
\beg{align*}&-\ff 1 2 \Ric_Z(\nn f,\nn f )(x)\\
 &\le \ff  1 2 K(x)|\nn f(x)|^2+ \limsup_{T\to 0} \ff 1 T \EE \Big\<\nn f(x), \nn P_T f(x) -\EE [U_0^x (U_T^x)^{-1} \nn f(X_T^x)]\Big\>.\end{align*} Since by  \eqref{GR2}, \eqref{QQ} and \eqref{LI}   we have
\beg{align*}&\big\<\nn f(x), \nn P_Tf(x)- \EE [U_0^x (U_T^x)^{-1} \nn f(X_T^x)]\big\>\\
&= -\int_0^T\big\<\nn f(x), U_0^x  \Ric_Z(U_r^x)(U_T^x)^{-1}\nn f(X_T^x)\big\>\d r =- T\Ric_Z(\nn f,\nn f)(x)+ {\rm o} (T)\end{align*}   for small $T>0$, this implies
$\Ric_Z(\nn f,\nn f)(x)\le K(x).$

On the other hand, to prove the desired bound on $\|\II\|$, we let $x\in \pp M$, $f\in C_0^\infty(M)$ with $\<\nn f, N\>(x)=0, |\nn f(x)|=1$ and $\Hess_f(x)=0$. By  \cite[Lemma 3.1.2]{W2},
$$\EE\e^{\ll l^x_{T\land\tau_1}}<\infty,\ \ \EE l_{T\land \tau_1}^x = \ff{2\ss T}{\ss\pi}+ {\rm O}(T^{3/2})$$ for all $\ll>0$ and small $T>0$.
Combining this with \eqref{MU}, \eqref{DE},  and \eqref{LO}, we obtain
\beq\label{LOC}\lim_{T\to 0}\ff {\EE\mu_{x,T}([0,T])} {\ss T} = \ff{2\si(x)}{\ss\pi},\ \ \  \lim_{T\to 0} \ff{[\EE \mu_{x,T}([0,T])]^2}{\ss T} = 0.\end{equation}
Then repeating the above argument with \eqref{II}  replacing \eqref{RIC}, we prove $$|\II(\nn f,\nn f)(x)|\le \si(x).$$ Indeed, by \eqref{II} and \eqref{LOC}, instead of \eqref{X1} we have
$$-\II(\nn f,\nn f)(x) \le \ff{\ss\pi}{2} \lim_{T\to\infty} \ff{|\nn P_T f|^p(x)- P_T|\nn f|^p(x)}{p\ss T} =\si(x),$$ while multiplying \eqref{X2} by $\ss T$ and letting $T\to\infty$ leads to
$$-\ff 1 {\ss\pi} \II(\nn f,\nn f)(x) \le \ff{\si(x)}{\ss\pi} -\ff 2 {\ss\pi} \II(\nn f,\nn f)(x).$$

{\bf (d)}  (5) $\Rightarrow$ (1). Let $F(\gg)= f(\gg_T)$. Then (5) implies
\beq\label{55} P_Tf^2(x)-  (P_Tf(x))^2 \le 2 \int_0^T \EE \big[(1+\mu_{x,T}([s,T]))^2 |\nn f(X_T^x)|^2\big]\d s.\end{equation}  For $f$ in \eqref{RIC}, combining this with \eqref{RIC} and \eqref{MUL} we obtain
\beg{align*} &\Ric_Z(\nn f, \nn f)(x) =\lim_{T\to 0} \ff 1 T \bigg(\ff{P_T f^2(x)- (P_Tf)^2(x)}{2T}-|\nn P_T f|^2\bigg)\\
&\le \lim_{T\to 0}\ff 1 T \bigg\{\ff 1 T\int_0^T \Big\{\EE\big[(1+\mu([s,T]))^2|\nn f(X_T^x)|^2\big] -|\nn P_T f(x)|^2\Big)\d s\bigg\}\\
&=  \lim_{T\to 0}\ff 1 T \bigg\{P_T|\nn f|^2(x) -|\nn P_Tf|^2(x) +  \ff {2|\nn f|^2(x)}  {T}\int_0^T (T-s)K(x) \d s \bigg\}\\
&= 2 \Ric_Z(\nn f,\nn f)(x) + K(x)|\nn f|^2(x).\end{align*}
This implies $\Ric_Z(\nn f,\nn f)(x)\ge -K(x)|\nn f(x)|^2$.  Next, for $f$ in \eqref{II}, combining \eqref{55} with \eqref{II} and \eqref{LOC},  we obtain
\beg{align*} & \II(\nn f, \nn f)(x) =\lim_{T\to 0} \ff {3\ss\pi} {8\ss T} \bigg(\ff{P_T f^2(x)- (P_Tf)^2(x)}{2T}-|\nn P_T f(x)|^2\bigg)\\
&\le \lim_{T\to 0}\ff {3\ss\pi} {8\ss T} \bigg\{\ff 1 T\int_0^T \Big\{\EE\big[(1+\mu([s,T]))^2|\nn f(X_T^x)|^2\big] -|\nn P_T f(x)|^2\Big)\d s\bigg\}\\
&=  \lim_{T\to 0}\ff  {3\ss\pi} {8\ss T} \bigg\{P_T|\nn f|^2(x) -|\nn P_Tf|^2(x) +  \ff {2|\nn f(x)|^2}{T}\int_0^T \ff{2\si(x)(\ss T-\ss s)}{\ss\pi}   \d s +{\rm o}\big(\ss T\big)\bigg\}
 \\ &= \ff 3 2 \II(\nn f,\nn f)(x)+ \ff 1 2 \si(x).  \end{align*} Hence, $\II(\nn f,\nn f)(x)\ge -\si(x)|\nn f(x)|^2.$

On the other hand, to prove the upper bound estimates, we take $F(\gg)= f(\gg_\vv)-\ff 1 2f(\gg_T)$ for $\vv\in (0,T).$  By \eqref{MG},
$$|\dot D_t F|= \Big|\nn f(X_\vv)- \ff 1 2 U_\vv^x(U_T^x)^{-1}\nn f(X_T^x)\Big|1_{[0,\vv)}(t)+ \ff 1 2 |\nn f(X_T^x)|1_{[\vv,T]}(t).$$
Then (5) implies
\beq\label{56}\beg{split} I_\vv &:= \EE\Big[f(X_\vv^x)-\ff 1 2 \EE(f(X_T^x)|\F_\vv)\Big]^2-\Big(P_\vv f(x)-\ff 1 2P_Tf(x)\Big)^2\\
&\le 2  \vv \EE\bigg\{(1+\mu_{x,T}([0,T]))\Big(\Big|\nn f(X_\vv^x)- \ff 1 2 U_\vv^x(U_T^x)^{-1}\nn   f(X_T^x)\Big|^2 \\
&\qquad\quad \  +\ff {\mu_{x,T}([0,T])|\nn f(X_T^x)|^2} 4\bigg\}    +c \vv^2 =: J_\vv,\ \ \vv\in (0,T)\end{split}\end{equation} for some constant $c>0$. Obviously,
\beq\label{57}\beg{split}  \lim_{\vv\to 0} \ff{J_\vv}\vv = \EE\bigg\{&(1+\mu_{x,T}([0,T]))\Big(\Big|\nn f(x) - \ff 1 2 U_0^x(U_T^x)^{-1}\nn   f(X_T^x)\Big|^2\\
&+\ff {\mu_{x,T}([0,T])}4 |\nn f|^2(X_T^x)\Big)\bigg\}.\end{split}\end{equation}

On the other hand, we have
\beq\label{IV} \beg{split} \ff{I_\vv}\vv =& \ff{P_\vv f^2-(P_\vv f)^2}\vv +\ff 1{4\vv} \EE\Big[\big\{\EE\big(f(X_T^x)|\F_\vv\big)\big\}^2 -(P_Tf)^2(x)\Big]\\
& + \ff{\EE[f(X_T^x)\{P_\vv f(x)-f(X_\vv^x)\}]}\vv.\end{split} \end{equation} Let $f\in C_0^\infty(M)$ satisfy the Neumann boundary condition, we have
\beq\label{IV2} \lim_{\vv\to 0} \ff{P_\vv f^2-(P_\vv f)^2}\vv= 2|\nn f|^2(x). \end{equation}
Next, \eqref{TTD} and  \eqref{MF} yield
\beq\label{*0} \EE(f(X_T^x)|\F_\vv)= P_Tf(x)+\ss 2 \int_0^\vv \big\<\EE\big(Q_{s,T}^x(U_T^x)^{-1}\nn f(X_T^x)\big|\F_s\big), \d W_s\big\>.\end{equation}
Then $$\EE[\EE(f(X_T^x)|\F_\vv)]^2= (P_Tf)^2 + 2 \int_0^\vv\EE |Q_{0,T}^x(U_T^x)^{-1}\nn f(X_T^x)|^2\d s.$$ This together with \eqref{GR2} leads to
\beq\label{IV3}\beg{split}  & \lim_{\vv\to 0}\ff 1{4\vv} \EE\Big[\big\{\EE\big(f(X_T^x)|\F_\vv\big)\big\}^2 -(P_Tf)^2(x)\Big]\\
 &=  \ff 1 2 \Big|\EE\big[Q_{0,T}^x(U_T^x)^{-1}\nn f(X_T^x)\big]\Big|^2 = \ff 1 2 |\nn P_T f(x)|^2.   \end{split}\end{equation}
 Finally, by It\^o's formula we have
\beg{align*} P_\vv f(x)-f(X_\vv^x)& = P_\vv f(x)-f(x) -\int_0^\vv Lf(X_s^x)\d s -\ss 2 \int_0^\vv \<\nn f(X_s^x), U_s^x\d W_s\>\\
&= {\rm o}(\vv)  -\ss 2 \int_0^\vv \<\nn f(X_s^x), U_s^x\d W_s\>.\end{align*} Combining this with \eqref{*0} and \eqref{GR2}, we arrive at
$$\lim_{\vv\to 0} \ff{\EE[f(X_T^x)\{P_\vv f(x)-f(X_\vv^x)\}]}\vv=- 2\<\nn f(x),\nn P_t f(x)\>.$$
Substituting this and \eqref{IV2}-\eqref{IV3} into \eqref{IV}, we obtain
$$\lim_{\vv\to 0}\ff{I_\vv}\vv= 2\Big|\nn f(x)- \ff 1 2 \nn P_Tf(x)\Big|^2.$$
Combining this with \eqref{56} and \eqref{57}, we prove the second inequality in (2) for $q=2$, which implies $\Ric_Z\le K$ and $\II\le \si$ as shown in step {\bf (c)}.

{\bf (e)} $(1)\Rightarrow (4).$ According to \eqref{MF},
\begin{equation}\label{eq2.27}G_t:= \EE(F^2|\F_t) =\mathbb{E}(F^2)+\sqrt{2}\int^t_0\big\<\EE(\tt D_s F^2|\F_s),\d W_s\big\>,\ \ t\in [0,T].\end{equation}
By It\^o's formula,
\begin{equation}\label{eq2.28}\aligned \d (G_t\log G_t)&=(1+\log G_t)\d G_t+\frac{|\EE(\tt D_s F^2|\F_s)|^2}{G_t}\d t\\
&\le  (1+\log G_t) \d G_t +4 \EE(|\tt D_sF|^2|\F_s)\d t.\endaligned\end{equation}
Then
\beq\label{LST} \EE[G_{t_1}\log G_{t_1}]-\EE[G_{t_0}\log G_{t_0}]
\le 4 \int_{t_0}^{t_1}  \EE  |\tt D_s F|^2\d s .\end{equation}
By (\ref{TTD})  we have
$$\aligned&\tt D_sF=\sum_{i=1}^N1_{\{s< t_i\}}Q^x_{s,t_i}(U_{t_i}^x)^{-1}\nabla_if\\
&=\sum_{i=1}^N1_{\{s< t_i\}}\bigg(I-\int^{t_i}_sQ^x_{s,t}\big\{ \Ric_V(U_t^x)\d t+
 \mathbb{I}_{U^x_t}\d l_t^x\big\}\bigg)\left(I-1_{\{X^x_{t_i}\in \partial M\}}P_{U^x_{t_i}}\right)(U_{t_i}^x)^{-1}\nabla_if\\
&=\dot{D}_0F-\int^T_sQ^x_{s,t} \big\{\Ric_Z(U_t^x)\d t+
 \II(U^x_t)\d l_t^x\big\}.\endaligned$$
Combining this   with (1), \eqref{QQ2} and \eqref{MUL}, and using the Schwarz inequality, we prove
\begin{equation}\label{eq2.25}|\tt D_sF|^2 \le (1+\mu_{x,T}([s,T])) \bigg(|\dot D_0F|^2 +\int_s^T |\dot D_sF|^2\mu_{x,T}(\d s)\bigg).\end{equation}
This together with \eqref{LST} implies the log-Sobolev inequality in (4).

 \end{proof}

\section{Extension of Theorem \ref{T1.1}}

In this section, we aim to drop the condition \eqref{DE} in Theorem \ref{T1.1} and allow  the (reflecting) diffusion process generated by $L$ to be explosive.
The idea is to  make a conformal change of metric such that the condition \eqref{DE} holds on the new Riemannian manifold.
Since  both $\Ric_Z$ and
$\II$ are local quantity, they doe  not change at $x$ if the new metric coincides with the original one around point $x$.

Let $(M,g)$ be a Riemannian manifold with boundary, and let $N$ be the inward pointing unit normal vector field of $\partial M$. Let $\phi\in
C_0^\infty(M)$ be non-negative with non-empty  $M_\phi:= \{\phi>0\}.$
  Then, $M_\phi$ is a complete
Riemannian manifold under the metric $g_\phi:= \phi^{-2}g$. Let $\nn^\phi, \DD^\phi, \Ric^\phi$ and $\II^\phi$ be the associated Laplacian, gradient, Ricci curvature and the second fundamental form of $\pp M_\phi$.
By  e.g. \cite[Theorem 1.159 d)]{E},
$$\nn^\phi_X Y= \nn_XY -\<X, \nn\log \phi\>Y - \<Y, \nabla\log\phi\>X +\<X,Y\>\nn\log \phi.$$
Moreover, according to \cite[Theorem 1.2.4]{W2} and the proof of \cite[Theorem 1.2.5]{W2},  we have
\beg{align*} &\Ric_\phi= \Ric +(d-2)\phi^{-1}\Hess_\phi +(\phi^{-1}\DD \phi -(d-3)|\nn \log
\phi|)g,\\
&\II^\phi =\phi^{-1} \II + (N\log\phi) g.\end{align*}
Noting that  $|X|=1$ if and only if $g_\phi(\phi X, \phi X)=1$, we obtain
$$\|\II_g\|_\infty =\sup_{X\in T\pp M_\phi, |X|=1} |\II_\phi(\phi X, \phi X)|<\infty,$$
and for $\Ric_{\phi Z}^\phi$ the curvature of $L^\phi:= \DD^\phi+ \phi Z$,
$$\|\Ric_{\phi Z}^\phi\|_\infty= \sup_{X\in TM_\phi, |X|=1} |\Ric^\phi(\phi X, \phi X)- g_\phi(\nn_{\phi X}(\phi Z), \phi X)|<\infty.$$
Therefore, Theorem \ref{T1.1} applies to $L^\phi$ on the manifold $M_\phi$. In particular, by taking $\phi$ such that $\phi=1$ around a point $x$, we have $\Ric_Z=\Ric^\phi$ and $\II=\II^\phi$ at point $x$, so that in this way we characterize these two quantities at $x$. To this end,  we will take $\phi= \ell(\rr_x)$, where $\rr_x$ is the Riemannian distance to $x$ and $\ell\in C_0^\infty(\R)$ is such that $0\le\ell\le 1$, $\ell(s)=1$ for $s\le r$ and $\ell (s)=0$ for $s\ge 2r$ for some constant $r>0$ with compact $B_{2r}(x):= \{\rr_x\le 2r\}$.

Obviously,  before exiting the ball $B_r(x)$  the   diffusion process generated by $L$ coincides with that generated by $L^\phi$. So, to use the original diffusion process in place of the new one,  we will take references functions which vanishes as soon as the diffusion exits this ball. To this end, we will make truncation of cylindrical functions in terms of the uniform distance
\begin{equation*}
\tilde{\rho_x}(\gamma):=\displaystyle\sup_{t\in[0,1]}\rho(\gamma(t),x).
\end{equation*}
To make the manifold $M_\phi$ complete,  let  $\dd: M\to (0,\infty)$ be a smooth function such that $B_R(x) $ is compact for any $R\le \dd_x$.  Consider the class of truncated cylindrical functions
\begin{equation}\label{e1}
\F C^\infty_{T,loc}:=\Big\{F\ell(\tilde{\rho_x}): F \in \F C^\infty_{T},\ x\in M, \ \ell \in C_0^{\infty}(\R), \ {\rm supp}\ell\subset [0,\dd_x)\Big\}.
\end{equation}

To define $\E_{t,T}^{K,\si}(\tt F,\tt F)$ for $\tt F=F\ell(\tilde{\rho_x})\in \F C^\infty_{T,loc}$, we take $\phi\in C_0^\infty(M)$ such that $0\le \phi\le 1,\
\phi=1$ for $\ell(\rr_x)>0$, and $\phi=0$ for $\rr_x\ge \dd_x$. Then $M_\phi$ is complete with bounded $\Ric_{\phi Z}^\phi$ and $\II^\phi$. Let $X_{[0,T]}^{x,\phi}$ be the (reflecting) diffusion process generated by $L^\phi$.  Similarly to the proof of \cite[Lemma 2.1]{CW} for the case without boundary,  we see that $|\dot D_s \tt F(X_{[0,T]}^{x,\phi})|$ is well defined and bounded for $s\in [0,T].$ Noting that $\tt F$ is supported on  $\{\ell(\tt\rr_x)>0\}\subset W_T(M^\phi)$ and  $X_{[0,T]}^{x,\phi}=X_{[0,T]}^x$ if $\ell(\tt\rr_x(X_{[0,T]}^{x,\phi}))>0$ (see \eqref{eq3.3} below), we conclude that $|\dot D_s \tt F(X_{[0,T]}^{x})|=|\dot D_s \tt F(X_{[0,T]}^{x,\phi})|$ is well defined and bounded in $s\in [0,T]$ as well, which does not depend on the choice of $\phi$. Again since $\tt F$ is supported on $\{\ell(\tt\rr_x)>0\}\subset W_T(M^\phi)$  and $M^\phi$ is relatively compact in $M$,    we have
$$\E^{K,\sigma}_{t,T}(\tt F,\tt F)
 := \EE\bigg\{\big(1+\mu_{x,T}([t,T])\big)
  \bigg(|\dot{D}_t\tt F(X_{[0,T]}^x)|^2
 +\int^{T}_{t}   |\dot{D}_s\tt F(X_{[0,T]}^x)|^2\mu_{x,T}(\d s)\bigg)\bigg\}<\infty.$$

\beg{thm}\label{T3.1}  Let $K\in C(M; [0,\infty))$ and $\si\in C(\pp M; [0,\infty))$. The following statements are equivalent each other:
\beg{enumerate}\item[$(1)$] For any $x\in M$ and $y\in\pp M$,
\beg{align*} &\|\Ric_Z\|(x):= \sup_{X\in T_x M, |X|=1} |\Ric(X,X)-\<\nn_XZ, X\>|(x)\le K(x),\\
&\|\II\|(y):= \sup_{Y\in T_y\pp M, |Y|=1}|\II(Y,Y)|(y) \le \si(y).\end{align*}
\item[$(2)$] For any  $t_0,t_1\in [0,T]$ with $t_1>t_0$, and any $ x\in M$,  the following log-Sobolev inequality holds:
$$\aligned&\EE\left[\EE\big(F^2(X_{[0,T]}^x)|\F_{t_1}\big) \log \EE(F^2(X_{[0,T]}^x)|\F_{t_1})\right]\\&-\EE\left[\EE\big(F^2(X_{[0,T]}^x)|\F_{t_0}\big) \log \EE(F^2(X_{[0,T]}^x)|\F_{t_0})\right]
 \le  4 \int_{t_0}^{t_1} \E_{s,T}^{K,\si}(F,F)\d s,\  \  F\in \F C^\infty_{T,loc}.\endaligned$$
 \item[$(3)$] For any  $t \in [0,T]$  and $ x\in M$,  the following Poincar\'e inequality holds:
$$   \EE \Big[\big\{\EE(F(X_{[0,T]}^x)|\F_{t})\big\}^2\Big]-   \Big\{\EE  \big[F (X_{[0,T]})\big] \Big\}^2
 \le 2\int_0^t \E_{s,T}^{K,\si}(F, F)\d s,\ \ F\in \F C^\infty_{T,loc}.$$   \end{enumerate}\end{thm}

\begin{proof} Since  $(2)\Rightarrow (3)$ is well known,   we only prove $(1)\Rightarrow (2)$  and $(3)\Rightarrow (1)$.

{\bf (a)}  (1) $\Rightarrow$ (2). Fix $x\in M$. For any $\tilde{F}:=F\ell(\tilde{\rho}_x)\in \F C^\infty_{T,loc}$, there exists $R\in (0,\dd_x)$ such that $\supp(\ell(\tt\rho_x))\subset B_R(x):=\{y\in M: \rho(x,y)\leq R\}$. Let $\phi_R \in C_0^{\infty}(M)$
 such that $\phi_R|_{B_R(x)}=1$ and $0\le \phi_R\le 1$.    We consider the following Riemannian  metric  on the manifold $M_R:=\{y \in M:\ \phi_R(y)>0\}$:
\begin{equation*}
g_R:=\phi_R^{-2} g.
\end{equation*}
As explained above that $(M_R,g_R)$ is a complete
Riemannian manifold with
\begin{equation}\label{eq3.1}
K_R:=\sup_{M_R}\|\Ric_Z^{R}\|_\infty <\infty,\quad \sigma_R:=\sup_{M_R}\|\mathbb{I}^{R}\|_\infty <\infty.
\end{equation}

We consider the SDE (\ref{eq1.1}) on $M$,
\begin{equation}\label{eq3.2}
\begin{cases}
&\d U_t^x=\sqrt{2}\, H_{U_t^x}(U_t^x)\circ\d W_t+H_{Z}(U_t^x)\d t+H_N(U_t^x)\d l_t^x,\\
& U_0=u_0.
\end{cases}
\end{equation}
Then $X_t:=\pi(U_t)$ is the (reflecting if $\pp M$ exists) diffusion process on $M$ generated by $L=\DD+Z$.

Similarly, let $\{H_{i,R}\}_{i=1}^n$ and $H_{\phi_RZ,R}$ be the orthonormal basis   of horizontal vector fields and horizontal lift of $\phi_RZ$ under the metric $g_R$. Since $g_R=g$ and $\phi_R=1$ on $B_R(x)$, for $u\in O(M_R) $ with $\pi u\in B_R(x)$ we have
  $H_{i,R}(u) =H_i(u)$ and $H_{\phi Z,R}(u)=  H_Z(u)$.  For  $W_t$ and $u_0$   in (\ref{eq3.2}),
we consider the following SDE on the manifold $M_R$:
\begin{equation*}
\begin{cases}
&\d U_{t,R}=\displaystyle\sum^n_{i=1}H_{i,R}(U_{t,R})\circ\d W_t^i+H_{\phi_RZ,R}(U_t^x)\d t+H_N(U_t^x)\d l_{R, t}^x,\\
& U_{0,R}=u_0.
\end{cases}
\end{equation*}
Then $X^{x,R}_{\cdot}:=\pi(U_{\cdot,R})$ is the (reflecting if $\pp M_R$ exists) diffusion process on $M_R$ generated by $L_R:=\DD_R+\phi_R Z$, where $\DD_R$ is the Laplacian on $M_R$.  Obviously,
\begin{equation}\label{eq3.3}U_{t,R}= U_t, \ l_{R, t}^x= l_{t}^x \ \ \text{for}  t\le  \tau_R:=
\inf\{t\ge 0: X_t \notin B_R(x)\}.
 \end{equation}
Denote by $\P^T_{R,x}$ the distribution of
the  process $X^{x,R}_{[0,T]}$. By \cite{W2} and (\ref{LST}), we have  the damped logarithmic Sobolev inequality holds
\begin{equation}\label{eq3.4}
\EE[G_{t_1}\log G_{t_1}]-\EE[G_{t_0}\log G_{t_0}]
\le 4 \tilde{\E}^{t_1,t_0}_R(G,G),\quad G\in \F C^\infty_{T},
\end{equation}
where $G_t:= \EE(G^2(X^{x,R}_{[0,T]})|\F_t)$ and
$$\tilde{\E}^{t_1,t_0}_R(H,G)=\int_{W^T_x(M_R)}\int^{t_1}_{t_0}\<\tilde{D}^R_sF,\tilde{D}^R_sG\>\d s\d \P^T_{R,x}.$$
According to \cite{W2},  the form $(\tilde{\E}^{t_1,t_0}_R,\F C^\infty_{T})$ is closable in $L^2(\P^T_{R,x})$. Let   $(\tilde{\E}^{t_1,t_0}_R,\D(\tilde{\E}^{t_1,t_0}_R))$ be its closure.
Let $\rho^R$ be the Riemannian distance on $M_R$ and
$$\tilde{\rho_x}^R(\gamma):=\sup_{t \in [0,1]}\rho^R(\gamma(t),x),\quad \gamma \in W^T_x(M_R).$$
We have $\tilde{\rho_x}^R(\gamma)=\tilde{\rho_x}(\gamma)$
for each $\gamma \in W^T_x(M_R)\subseteq W^T_x(M)$ satisfying
$\rho_x^R(\gamma)\le R$.
Then  \cite[Lemma 2.1]{CW} implies that  $\ell(\tilde{\rho}_x)$ is in $ \D(\tilde{\E}_{\P^T_{R,x}})$, and so is $\tilde{F}:=F\ell(\tilde{\rho}_x)$. Combining this with (\ref{eq3.3}) and (\ref{eq3.4}), we get
\begin{equation}\label{eq3.6}\aligned &\EE\bigg[\EE\big(\tt F^2(X_{[0,T]}^x)|\F_{t_1}\big) \log\EE(\tt F^2(X_{[0,T]}^x)|\F_{t_1})\bigg]\\
&\quad -\EE\bigg[\EE\big(\tt F^2(X_{[0,T]}^x)|\F_{t_0}\big) \log\EE(\tt F^2(X_{[0,T]}^x)|\F_{t_0})\bigg]\\
&=\EE\bigg[\EE\big(\tt F^2(X_{[0,T]}^{x,R})|\F_{t_1}\big) \log\EE(\tt F^2(X_{[0,T]}^{x,R})|\F_{t_1})\bigg]\\
&\quad -\EE\bigg[\EE\big(\tt F^2(X_{[0,T]}^{x,R})|\F_{t_0}\big) \log\EE(\tt F^2(X_{[0,T]}^{x,R})|\F_{t_0})\bigg]
\\&\leq4 \int_{W^T_x(M_R)}\int^{t_1}_{t_0}\<\tilde{D}^R_s\tt F,\tilde{D}^R_s\tt F\>\d s\d \P^T_{R,x}=4 \int_{W^T_x(M)}\int^{t_1}_{t_0}\<\tilde{D}_s\tt F,\tilde{D}_s\tt F\>\d s\d \P^T_{x}.\endaligned\end{equation}
Combining this with  (\ref{eq2.25}), we prove (2).

{\bf (a)}  (3) $\Rightarrow$ (1). We first prove the lower bound estimates.
When   $x\in M\setminus \pp M$, there exists $r\in (0,\ff 1 2 \dd_x)$ such that $B_{2r}(x)\subset M\setminus \pp M$.
Let $\Phi=\ell(\tilde{\rho_x})$, where $\ell\in C_0^\infty(\R)$ such that $0\le \ell\le 1$, $\ell(s)=1$ for $s\le r$ and $\ell(s)=0$ for $s\ge 2r.$ Let $\tau_s=\inf\{t\ge 0: \rr(x,X_t^x)\ge s\}$ for $s>0$.
Consider  $\tilde{F}(\gg)= (\Phi F)(\gg)= \Phi(\gg) f(\gg_T)$ for $f$ in \eqref{RIC}.  Then (3) and \eqref{LO} imply
\beq\label{LWW}\beg{split}
&\EE \Big[(F\Phi)^2(X_{[0,T]}^x)\Big]-   \Big\{\EE  \big[(F\Phi) (X_{[0,T]})\big] \Big\}^2 \le 2\int_0^T \E_{t,T}^{K,\si}(\tt F, \tt F)\d t\\
&=2\int_0^T\EE\bigg\{\big(1+\mu_{x,T}([t,T])\big)\bigg(|\dot{D}_t\tt F(X_{[0,T]}^x)|^2 +\int^{T}_{t}   |\dot{D}_s\tt F(X_{[0,T]}^x)|^2\mu_{x,T}(\d s)\bigg)\bigg\} \d t\\
&\le 2\int_0^T\EE\Big[1_{\{\tau_{2r}>T\}} \big(1+\mu_{x,T}([t,T])\big)^2|\nabla f(X_T^x)|^2 \Big] \d t+C\P(\tau_r\leq T)\\
&=2\int_0^T\EE\Big[1_{\{\tau_{2r}>T\}} \big(1+\mu_{x,T}([t,T])\big)^2|\nabla f(X_T^x)|^2 \Big] \d t+{\rm o}(T^3), \end{split}\end{equation}
where $C>0$ is a constant depending on $f$ and $\Phi$.
On the other hand, by \eqref{RIC} and   \eqref{LO}, we have
\beg{align*}
&\lim_{T\to 0} \ff 1 T \bigg(\ff{\EE[F^2\Phi^2(X_{[0,T]}^x)]- \left\{\EE  \big[F\Phi(X_{[0,T]})\big] \right\}^2}{2T}-|\nn P_T f|^2\bigg)\\
&=\lim_{T\to 0} \ff 1 T \bigg(\ff{P_T f^2(x)- (P_Tf)^2(x)}{2T}-|\nn P_T f|^2\bigg)\\
&=\Ric_Z(\nn f, \nn f)(x). \end{align*}
Since $l_s^x=0$ for $s\le \tau_{2r}$, these two estimates together with \eqref{LO} and \eqref{MU} lead to
\beg{align*} &\Ric_Z(\nn f, \nn f)(x) =\lim_{T\to 0} \ff 1 T \bigg(\ff{\EE[(F\Phi)^2(X_{[0,T]}^x)]- \left\{\EE  \big[(F\Phi)(X_{[0,T]})\big] \right\}^2}{2T}-|\nn P_T f|^2\bigg)\\
&\le \lim_{T\to 0}\ff 1 T \bigg\{\ff 1 T\int_0^T \Big\{\EE\big[1_{\{\tau_{2r}>T\}} (1+\mu([s,T]))^2|\nn f(X_T^x)|^2\big] -|\nn P_T f(x)|^2\Big)\d s\bigg\}\\
&\le\lim_{T\to 0}\bigg(\ff {P_T|\nn f|^2(x) -|\nn P_Tf|^2(x)}{T} +\ff{\int_0^T\EE  \{1_{\{\tau_{2r}>T\}} [(1+\mu([s,T]))^2-1 ] |\nn f(X_T^x)|^2 \}\d s}{T^2}\bigg)\\
&= 2 \Ric_Z(\nn f,\nn f)(x) + K(x)|\nn f|^2(x). \end{align*}
 Therefore,   $\Ric_Z(\nn f,\nn f)(x)\ge -K(x)|\nn f(x)|^2$.

Next, let $x\in \pp M$. For $f$ in \eqref{II},  by \eqref{LO} we have
\beq\label{eq3.7}\aligned
&\lim_{T\to 0} \ff {3\ss\pi} {8\ss T}\bigg(\ff{\EE[(F\Phi)^2(X_{[0,T]}^x)]- \left\{\EE  \big[(F\Phi)(X_{[0,T]})\big] \right\}^2}{2T}-|\nn P_T f|^2\bigg)\\
&=\lim_{T\to 0} \ff {3\ss\pi} {8\ss T}\bigg(\ff{P_T f^2(x)- (P_Tf)^2(x)}{2T}-|\nn P_T f|^2\bigg)\\
&=\II(\nn f, \nn f)(x).\endaligned\end{equation}
Combining this with \eqref{LWW} and \eqref{LOC},  we obtain
\beg{align*} & \II(\nn f, \nn f)(x) =\lim_{T\to 0} \ff {3\ss\pi} {8\ss T} \bigg(\ff{\EE[(F\Phi)^2(X_{[0,T]}^x)]- \left\{\EE  \big[(F\Phi)(X_{[0,T]})\big] \right\}^2}{2T}-|\nn P_T f(x)|^2\bigg)\\
&\leq \lim_{T\to 0} \ff {3\ss\pi} {8\ss T} \bigg(\int_0^T\ff{\EE\big\{1_{\{\tau_{2r}>T\}}\big(1+\mu_{x,T}([t,T])\big)^2 |\nabla F(X_T^x)|^2\big\}}{T}\d t-|\nn P_T f(x)|^2\bigg)
\\&=\lim_{T\to 0}\ff  {3\ss\pi} {8\ss T} \bigg\{P_T|\nn f|^2(x) -|\nn P_Tf|^2(x) +  \ff {2|\nn f(x)|^2}{T}\int_0^T \ff{2\si(x)(\ss T-\ss s)}{\ss\pi}   \d s  \bigg\}
 \\ &= \ff 3 2 \II(\nn f,\nn f)(x)+ \ff 1 2 \si(x).  \end{align*}
Therefore, $\II(\nn f,\nn f)(x)\ge -\si(x)|\nn f(x)|^2.$

To prove the upper bound estimates, we take $F(\gg)= f(\gg_\vv)-\ff 1 2f(\gg_T)$ for $\vv\in (0,T).$  By \eqref{MG},
$$|\dot D_t F|= \Big|\nn f(X_\vv)- \ff 1 2 U_\vv^x(U_T^x)^{-1}\nn f(X_T^x)\Big|1_{[0,\vv)}(t)+ \ff 1 2 |\nn f(X_T^x)|1_{[\vv,T]}(t).$$
Moreover, by (3) and \eqref{LO}, we may find a constant $C>0$ depending on $f$ and $\Phi$ such that for any $\vv,T\in (0,1)$,
\beq\label{eq3.8}\aligned
I_\vv&:= \EE\Big[ \EE\Big(\Phi(X_{[0,T]}^x)f(X_\vv^x)-\ff 1 2\Phi(X_{[0,T]}^x)f(X_T^x)\Big|\F_\vv\Big)\Big]^2\\
&~~~~~~~~~~~~~~~~~-\Big[ \EE\Big(\Phi(X_{[0,T]}^x)f(X_\vv^x)-\ff 1 2\Phi(X_{[0,T]}^x)f(X_T^x)\Big)\Big]^2\\&
\le 2\int^\varepsilon_0\EE\bigg\{\big(1+\mu_{x,T}([t,T])\big)
   |\Phi(X_{[0,T]}^x)\dot D_t F|^2
\\&~~~~~~~~~~~~~~~~ +\int^{T}_{t}   |\Phi(X_{[0,T]}^x)\dot D_s F |^2\mu_{x,T}(\d s)\bigg)\bigg\}\d t+ C\vv T^4.
\endaligned\end{equation}
Then
\beq\label{eq3.9}\beg{split}  \limsup_{\vv\to 0} \ff{I_\vv}\vv \le  \EE\bigg\{&\Phi(X_{[0,T]}^x)(1+\mu_{x,T}([0,T]))\Big(\Big|\nn f(x) - \ff 1 2 U_0^x(U_T^x)^{-1}\nn   f(X_T^x)\Big|^2\\
&+\ff {\Phi(X_{[0,T]}^x)\mu_{x,T}([0,T])}4 |\nn f|^2(X_T^x)\Big)\bigg\}+ {\rm o}(T^3) \end{split}\end{equation} for small $T>0$.
On the other hand, according to (d) of proof in Theorem \ref{T1.1}, we have
\beq\label{eq3.10} \beg{split} \ff{I_\vv}\vv =& \ff{P_\vv f^2-(P_\vv f)^2}\vv +\ff 1{4\vv} \EE\Big[\big\{\EE\big(f(X_T^x)|\F_\vv\big)\big\}^2 -(P_Tf)^2(x)\Big]\\
& + \ff{\EE[f(X_T^x)\{P_\vv f(x)-f(X_\vv^x)\}]}\vv+o(T^3\\
&= 2\Big|\nn f(x)- \ff 1 2 \nn P_Tf(x)\Big|^2+o(T^3). \end{split}\end{equation}
Combining this with  (\ref{eq3.9}),     we arrive at
\beq\label{eq3.11}\beg{split}
&2\Big|\nn f(x)- \ff 1 2 \nn P_Tf(x)\Big|^2\\
&\leq\EE\bigg\{\Phi(X_{[0,T]}^x)(1+\mu_{x,T}([0,T]))\Big(\Big|\nn f(x) - \ff 1 2 U_0^x(U_T^x)^{-1}\nn   f(X_T^x)\Big|^2\\
&+\ff {\Phi(X_{[0,T]}^x)\mu_{x,T}([0,T])}4 |\nn f|^2(X_T^x)\Big)\bigg\}+o(T^3)\end{split}
 \end{equation}
With this estimate, we may repeat  the last part in the proof of (2) $\Rightarrow$ (1) of Theorem \ref{T1.1} to derive the desired upper bound estimates on
$\Ric_Z$ and $\II$ at point $x$.
\end{proof}

\beg{thebibliography}{99}

\leftskip=-2mm
\parskip=-1mm

\bibitem{A} S. Aida, \emph{Logarithmic derivatives of heat kernels and logarithmic Sobolev inequalities
with unbounded diffusion coefficients on loop spaces,} J. Funct.
Anal. 174(2000), 430--477.

\bibitem{ATW09}  M. Arnaudon, A. Thalmaier, F.-Y. Wang, \emph{Gradient estimates and Harnack inequalities on non-compact Riemannian manifolds,} Stoch. Proc. Appl. 119(2009), 3653--3670.



\bibitem{CHL} B. Capitaine,  E. P. Hsu,  M. Ledoux, \emph{Martingale representation and a simple proof of logarithmic
Sobolev inequalities on path spaces,} Elect. Comm. Probab. 2(1997),
71--81.

\bibitem{CW} X. Chen, B. Wu, \emph{Functional inequality on path space
over a non-compact Riemannian manifold,}  J.
Funct. Anal. 266(2014), 6753-6779.

\bibitem{D} B. Driver,  \emph{A Cameron-Martin type quasi-invariant theorem
for Brownian motion on a compact Riemannian manifold},  J. Funct.
Anal. 110(1992), 272--376.


\bibitem{E} A. L. Besse, \emph{Einstein Manifolds,} Springer,
Berlin, 1987.

\bibitem{F} S. Fang, \emph{In\'egalit\'e du type de Poincar\'e sur l'espace
des chemins riemanniens}, C.R. Acad. Sci. Paris, 318 (1994),
257-260.




\bibitem{FWU} S. Z. Fang,  B. Wu, \emph{Remarks on spectral gaps on
the Riemannian path  space,}  {\it arXiv:1508.07657}.

\bibitem{HN} R. Haslhofer, A. Naber, \emph{Ricci curvature and Bochner formulas for martingales,} {\it arXiv:1608.04371}.

\bibitem{Hsu0} E. P. Hsu, \emph{Logarithmic Sobolev inequalities on path
spaces over Riemannian manifolds,} Comm. Math. Phys.  189(1997),
9--16.

\bibitem{H1} E. P. Hsu, \emph{Multiplicative functional for the heat equation on manifolds with boundary,} Mich. Math. J. 50(2002),351--367.



\bibitem{N} A. Naber, Characterizations of bounded Ricci curvature on smooth and nonsmooth spaces, {\it arXiv: 1306.6512v4}.

\bibitem{TW} A. Thalmaier and F.-Y. Wang, \emph{Gradient estimates for
harmonic functions on regular domains in Riemannian manifolds,} J. Funct. Anal. 155:1(1998),109--124.

\bibitem{RS} M. R\"{o}ckner and B. Schmuland, \emph{Tightness of general $C_{1,p}$ capacities on Banach
space,} J. Funct. Anal. 108(1992), 1--12.

\bibitem{W1} F.- Y. Wang, \emph{Weak poincar\'{e} Inequalities on path
spaces,} Int. Math. Res. Not. 2004(2004), 90--108.

\bibitem{W09} F.-Y. Wang, \emph{Second fundamental form and gradient of Neumann semigroups,}  J. Funct. Anal. 256(2009), 3461--3469.
\bibitem{W11}  F.-Y. Wang, \emph{Analysis on path spaces over Riemannian manifolds with boundary,} Comm. Math. Sci. 9(2011),1203--1212.

\bibitem{W2} F.- Y. Wang, \emph{Analysis for diffusion processes on Riemannian manifolds,} World Scientific, 2014.


\end{thebibliography}
\end{document}